\documentclass[11pt]{amsart}
\usepackage{amssymb}
\usepackage{enumitem}
\usepackage{graphicx}
\usepackage{xcolor}
\makeatletter
\@namedef{subjclassname@2020}{%
	\textup{2020} Mathematics Subject Classification}
\makeatother

\usepackage[T1]{fontenc}

\newtheorem{theorem}{Theorem}[section]
\newtheorem{corollary}[theorem]{Corollary}
\newtheorem{lemma}[theorem]{Lemma}
\newtheorem{proposition}[theorem]{Proposition}

\theoremstyle{definition}
\newtheorem{definition}[theorem]{Definition}
\newtheorem{remark}[theorem]{Remark}
\newtheorem{example}[theorem]{Example}

\numberwithin{equation}{section}

\begin{document}

	\baselineskip=17pt
 	\author[X.Q. Qiang]{Xiangqi Qiang}
	\address{ School of Science, Jiangsu University of Science and Technology,
		Zhenjiang 212100, China}
	\email{xq.qiang@just.edu.cn}

	\author[R. Liu]{Rui Liu}
\address{School of Mathematical Science, Yangzhou University,
	Yangzhou 225002, China}
\email{liurui20240048@163.com}	
 
	\author[C. J. Hou]{Chengjun Hou}
\address{School of Mathematical Science, Yangzhou University,
	Yangzhou 225002, China}
\email{cjhou@yzu.edu.cn}

		\thanks{Corresponding author: Chengjun Hou}

	\title[Invariance of Property (T) for \'Etale Groupoids]{Invariance of Property (T) for \'Etale Groupoids } 	
	
	\begin{abstract}
We study topological Property (T) for locally compact Hausdorff étale groupoids under pullbacks, amplifications, and Morita equivalencen. We show that Kazhdan pairs for the regular representation family pull back along  Haar system-preserving homomorphisms, and for general representation families we identify a constant-vector compatibility condition as the additional hypothesis needed. We also show that Property (T) for the regular representation family passes from a groupoid to its transformation groupoid with the same Kazhdan constant; under a natural restriction condition on representation families, the same holds for general families. Finally, we prove that Property (T) is invariant under ampliation for ample groupoids with compact unit spaces, both for all representations and for the regular representation family. Hence topological Property (T) is Morita invariant for ample Hausdorff étale groupoids with compact unit spaces, and likewise for the regular representation family.
	\end{abstract}
	
	
	
	\maketitle
	
	\normalsize

	\baselineskip=17pt
	\section{Introduction}

Kazhdan’s Property (T) is a rigidity property for locally compact groups, introduced by Kazhdan \cite{Kazhdan} in the mid-1960s to show that a large class of lattices are finitely generated. Since then, Property (T) has become a basic notion in domains as diverse as representation theory, group theory, differential geometry, ergodic theory, operator algebras, and the theory of algorithms; for further details, see \cite{Bekka} and the references therein. Quantitative versions such as Kazhdan constants have been developed for discrete groups \cite{Zuk}, and the notion has been extended to continuous homomorphisms between topological groups \cite{Leung}.

Measured analogues of Property (T) were developed by Zimmer \cite{Zimmer} for equivalence relations and by Anantharaman-Delaroche \cite{Anantharaman-Delaroche} in a more general setting.
For applications to groupoid $C^*$-algebras, however, a topological version is needed. Such a version for étale groupoids was introduced by Dell'Aiera and Willett \cite{DellAiera Willett}, simultaneously generalizing Kazhdan's Property (T) for groups and geometric Property (T) for coarse spaces.    One main goal of this notion  is to construct Kazhdan projections in maximal and reduced groupoid $C^*$-algebras, with applications to exactness, $K$-exactness, and the Baum–Connes conjecture.
 Dell'Aiera and Willett also proved that an HLS groupoid $G_{\mathcal{K}}$ with parent group $\Gamma$ has Property (T) if and only if $\Gamma$ does, and that, for a discrete group $\Gamma$ acting on a compact space $X$,  if $X$ admits an invariant probability measure, then the transformation groupoid $X\rtimes \Gamma$ has Property (T) if and only if $\Gamma$ does.  These results motivate our study of transformation groupoids beyond the invariant-probability-measure setting and, more generally, of pullback along morphisms of étale groupoids.

 In this paper, we study topological Property (T) for locally compact Hausdorff étale groupoids with compact unit spaces. In the group setting, Property (T) is inherited by closed subgroups of finite covolume (in particular, by lattices), passes to finite extensions, and is invariant under commensurability. These facts lead us to consider the analogous operations for groupoids: pullback along morphisms, amplification along surjective local homeomorphisms, and Morita equivalence. We develop a systematic theory of Property (T) under these operations, with transformation groupoids treated as an application of the pullback theorem. We prove an amplification theorem for ample groupoids with compact unit spaces and deduce Morita invariance of this notion of Property (T). This partially answers a question of Dell’Aiera and Willett \cite[Section 8, Question (iv)]{DellAiera Willett}.


Our first main result concerns the pullback of Property (T) along a   Haar system-preserving homomorphism $\phi:\mathcal{G}\to \mathcal{H}$ of \'etale groupoids with compact unit spaces. In the regular case, the pullback morphism not only intertwines the regular representations but also identifies their constant-vector spaces, so Kazhdan pairs for the regular representation family of $\mathcal{H}$ pull back to Kazhdan pairs for the regular representation family of $\mathcal{G}$,  without requiring surjectivity;  in particular, Property (T) for the regular representation family passes from $\mathcal{H}$ to $\mathcal{G}$.  For general representation families, constant vectors for $\pi$ are automatically constant for the   representation $\pi\circ\phi^{*}$, where $\phi^{*}$ is the dual map induced by $\phi$, but not conversely; this obstruction is precisely what the constant-vector compatibility condition is designed to overcome. More precisely, if $\mathcal{F}_{\mathcal{G}}$ and $\mathcal{F}_{\mathcal{H}}$ are families of representations of $C_c(\mathcal{G})$ and $C_c(\mathcal{H})$, respectively, and $(\mathcal{F}_{\mathcal{G}},\mathcal{F}_{\mathcal{H}})$ is constant-vector compatible with $\phi$, then Kazhdan pairs for $\mathcal{F}_{\mathcal{H}}$ pull back to Kazhdan pairs for $\mathcal{F}_{\mathcal{G}}$, and hence Property (T) passes from the  $(\mathcal{H},\mathcal{F}_{\mathcal{H}})$ to $(\mathcal{G},\mathcal{F}_{\mathcal{G}})$.   

%
 
As an application of the pullback theorem, we consider transformation
groupoids. Let $\mathcal{G}$ be a locally compact Hausdorff \'etale groupoid
with compact unit space, and let $X$ be a compact Hausdorff left
$\mathcal{G}$-space with (surjective) anchor map. The projection
$q:\mathcal{G}\ltimes X\to\mathcal{G}$, $q(g,x)=g$, is a  
Haar system-preserving homomorphism, so the pullback theorem applies. In the
regular case, Property (T) passes from $\mathcal{G}$ to
$\mathcal{G}\ltimes X$, and Kazhdan pairs pull back with the same Kazhdan constant.
More generally, if every representation in a family $\mathcal{E}$ of
$C_c(\mathcal{G}\ltimes X)$ restricts, up to unitary equivalence, to a
representation in a family $\mathcal{F}$ of $C_c(\mathcal{G})$, then
Property (T) passes from $(\mathcal{G},\mathcal{F})$ to
$(\mathcal{G}\ltimes X,\mathcal{E})$.  
Finally, we construct a counterexample  showing that the converse of the regular pullback theorem fails: Property (T) need not descend along surjective morphisms, even for full Property (T).

Our second main result is an amplification theorem. Let $\mathcal{G}$ be an ample \'{e}tale groupoid with compact unit space,
and let $\mathcal{G}^{\varphi}$ be its ampliation along a surjective local
homeomorphism $\varphi:X\to\mathcal{G}^{(0)}$ with $X$ compact Hausdorff.
Then $\mathcal{G}$ has Property (T) if and only if $\mathcal{G}^{\varphi}$
does, and the same equivalence holds for $(\mathcal{G},\mathcal{R}_{\mathcal{G}})$
and $(\mathcal{G}^{\varphi},\mathcal{R}_{\mathcal{G}^{\varphi}})$. This implies our main result, which asserts that 
topological Property (T) is invariant under Morita equivalence, and likewise for Property (T) of the regular representation family. Moreover, by the Farsi–Kumjian–Pask–Sims equivalence theorem, similarity, Renault equivalence, stable isomorphism, and weak Kakutani equivalence coincide with Morita equivalence; hence Property (T) is invariant under each of these equivalences as well.

%
  


 The paper is organized as follows. Section 2 recalls the relevant notions for groupoids and Property (T). Section 3 studies the pullback of Property (T) along Haar system-preserving morphisms, introduces the constant-vector compatibility condition, and proves the corresponding pullback theorem, while Section 4 treats transformation groupoids as an application of the pullback theorem. Section 5 establishes the amplification theorem and deduces Morita invariance from the definition of Morita equivalence.

	\section{Preliminaries}

	For a topological groupoid $\mathcal{G}$, let $\mathcal{G}^{(0)}$   be  the unit space. The range and source maps $r,s$ from $\mathcal{G}$ onto $\mathcal{G}^{(0)}$  are defined by $r(g)=gg^{-1}$ and $s(g)=g^{-1}g$, respectively. If $r$ and $s$ are local homeomorphisms, then $\mathcal{G}$ is called to be \emph{\'{e}tale}. An open subset $U \subseteq  \mathcal{G}$ is said to be an open bisection if both $r|_{U}$ and $s|_{U}$ are homeomorphisms onto their ranges. A groupoid $\mathcal{G}$  is ample if it is \'etale and $\mathcal{G} ^{(0)}$ is zero-dimensional; equivalently,
	$\mathcal{G}$  is ample if it has a basis of compact open bisections. We refer to \cite{Farsi,Renault,Sims} for more details on topological groupoids and their $C^*$-algebras.

	Throughout, $\mathcal{G}	$ will denote a locally compact Hausdorff \'etale groupoid whose unit space $\mathcal{G}^{(0)}$ is compact. For $u,v\in \mathcal{G}^{(0)}$, we write
	$
	\mathcal{G}	_u=s^{-1}(u)$ and $ \mathcal{G}	^u=r^{-1}(u)
	$	for the source and range  fibres over $u\in \mathcal{G}	^{(0)}$, respectively. Since $\mathcal{G}	$ is \'etale, these fibres are discrete, and then the counting measures on the source fibres form a Haar system.

For each $u\in \mathcal{G}^{(0)}$, the regular representation $	\pi_u:C_c(\mathcal{G})\longrightarrow B(\ell^2(\mathcal{G}_u))$
	is defined by
	\[
	\pi_u(f)\delta_\gamma
	=
	\sum_{\eta\in \mathcal{G}_{r(\gamma)}} f(\eta)\delta_{\eta\gamma},
	\qquad \gamma\in \mathcal{G}_u.
	\]
	We write $
	\mathcal{R}_\mathcal{G}=\{\pi_u: u\in \mathcal{G}^{(0)}\}$.
	
	For $f\in C_c(\mathcal{G}	)$, define
	\[	\Psi(f)(u)=\sum_{\gamma\in \mathcal{G}	^u} f(\gamma),\qquad u\in \mathcal{G}^{(0)}.	\]
	Then  $\Psi$ is a continuous linear from $C_c(\mathcal{G})$ onto $C(\mathcal{G}^{(0)})$. This map is also denoted by $\Psi_\mathcal{G}$ for stressing the groupoid $\mathcal{G}$.

	Let $(  H,\pi)$ be a  $*$-representation of $C_c(\mathcal{G}	)$. Following \cite{DellAiera Willett}, a vector $\xi\in  H$ is \emph{constant} or \emph{invariant} if
	\[
	\pi(f)\xi=\pi(\Psi	(f))\xi
	\]
	for all $f\in C_c(\mathcal{G}	)$. We denote by $  H^\pi$ the closed subspace of constant vectors, and set
	$
	  H_\pi=(  H^\pi)^\perp$.
	
	\begin{definition}	Let $\mathcal F$ be a family of $*$-representations of $C_c(\mathcal{G})$. A pair $(K,c)$, where $ K \subseteq \mathcal{G}$ is compact and $c>0$, is called a \emph{Kazhdan pair} of  $\mathcal{G}$ for $\mathcal F$ if, for every $(  H,\pi)\in\mathcal F$ and every $\xi\in H_\pi$, there exists $f\in C_c(\mathcal{G})$ with  support in $K$ and
		$\|f\|_I\le 1$ such that  $$
		\|\pi(f)\xi-\pi(\Psi(f))\xi\|\ge c\|\xi\|.$$

		The groupoid $\mathcal{G}$ is said to have \emph{topological Property {\rm(T)} with respect to $\mathcal F$} ( or \emph{$(\mathcal{G},\mathcal F)$ has    Property {\rm(T)}}) if it admits such a pair exists for $\mathcal F$. In particular,
  when  $\mathcal F$ is taken respectively to be the family of all $*$-representations of $C_c(\mathcal{G})$, the family $\mathcal{R}_\mathcal{G}$ of all regular $*$-representations of $C_c(\mathcal{G})$ or that of all $*$-representations of $C_c(\mathcal{G})$ which have on nonzero invariant vectors, the associated \emph{Property {\rm(T)}} is called \emph{full topological Property {\rm(T)}}, \emph{regular Property {\rm(T)}} or \emph{{\rm(}topological{\rm )} Property {\rm(T$_1$)}}, respectively.
\end{definition}
	
	From the definitions, we have a observation that full topological Property (T) and regular Property (T) are invariant under isomorphisms of \'etale groupoids. Example 5.2 in \cite{DellAiera Willett}   shows that Property (T) is strictly stronger than
Property (T$_1$).  The following example gives that Property (T) is strictly stronger than
regular Property (T).

\begin{example}


		
		Let $\Gamma=\mathrm{SL}_2(\mathbb Z)$ and let $\mathcal{K}=\{K_n\}$ be a nested sequence of principal congruence subgroups with trivial intersection, for instance $
		K_n=\ker\bigl(\mathrm{SL}_2(\mathbb Z)\to \mathrm{SL}_2(\mathbb Z/2^n\mathbb Z)\bigr)$.
		By Selberg's spectral-gap theorem, $\Gamma$ has Property $(\tau)$ with respect to the congruence subgroups, whereas $\Gamma$ does not have Kazhdan's Property (T).

		Let $\mathcal{G}_\mathcal{K}$ be the HLS groupoid associated to the approximating sequence $\mathcal{K}:=\{K_n\}$. By \cite[Proposition 4.15]{DellAiera Willett}, Property $(\tau)$ of the parent group $\Gamma$ implies that $\mathcal{G}_\mathcal{K}$ has Property (T) with respect to the family of representations that extend to $C^*_r(\mathcal{G}_\mathcal{K})$. Every regular representation belongs to this family, so $\mathcal{G}_\mathcal{K}$ has regular Property (T).	
		
		On the other hand,  \cite[Proposition 4.11]{DellAiera Willett} states that an HLS groupoid has full topological Property (T) if and only if its parent group has Kazhdan's Property (T). Since $\mathrm{SL}_2(\mathbb Z)$ does not have Property (T), $\mathcal{G}_\mathcal{K}$ does not have full topological Property (T).


\end{example}

	\section{Pullback of Property (T)}

	In this section, we study the pullback of Property (T) along  
	Haar system-preserving morphisms of \'etale groupoids. We first treat the regular case and then extend the discussion to general representation
	families, where an additional constant-vector compatibility condition is required.

	We now recall the following definitions and results from \cite[Definition 2.3, Remark 2.5(3), Propositions 2.6 and 2.7]{AustinMitra}.
Let $\mathcal{G}$ and $\mathcal{H}$ be two \'etale groupoids. A continuous groupoid homomorphism $\phi:\mathcal{ G}\to \mathcal{H}$ is called to be \emph{Haar system-preserving} if it is
	proper and  pushes the Haar measure on each source fibre on $\mathcal{G}$ forward to the corresponding Haar measure on $\mathcal{H}$. In our \'etale case, this is equivalent to fibrewise bijectivity: for each $u\in \mathcal{G}^{(0)}$,
	$
	\phi_u:\mathcal{G}_u\longrightarrow \mathcal{H}_{\phi(u)}$ is a bijection, and similarly for range fibres. Let $\phi^*$ be the dual map defined by $\phi^*(f)=f\circ\phi$ for $f\in C_c(\mathcal{H})$.

\begin{lemma}\label{le0} Let $\phi:\mathcal{ G}\to \mathcal{H}$ be a Haar system-preserving homomorphism. Then
\begin{enumerate}
\item[\rm(1)] $\phi^*$ is a $\ast$-homomorphism from $C_c(\mathcal{H})$ into $C_c(\mathcal{G})$, and is $I$-norm decreasing. In particular, $\phi^*(f)\in C(\mathcal{G}^{(0)})$ for $f\in C(\mathcal{H}^{(0)})$.

	 Moreover, $\phi^*$ extends to reduced completions; if $\phi$ is surjective, the induced map on reduced groupoid $C^*$-algebras is isometric.
\item[\rm(2)] For a representation $\pi: C_c(\mathcal{G})\to B(H)$, $H^{\pi}\subseteq H^{\pi\circ\phi^*}$ and $\Psi_{\mathcal{G}}(\phi^*(f))=\phi^*(\Psi_{\mathcal{H}}(f))$ for $f\in C_c(\mathcal{H})$.
\end{enumerate}	
\end{lemma}
\begin{proof} For (1), we refer to \cite{BonickeDellAieraGabeWillett}. For (2), we use the definition.
\end{proof}

In this section, we fix a Haar system-preserving homomorphism $\phi:\mathcal{ G}\to \mathcal{H}$.

	\subsection{Regular representations}

The following lemma comes from \cite{DellAiera Willett}.
	\begin{lemma}\label{Le1}
		Let $u\in \mathcal{G}^{(0)}$. The space of constant vectors of $\pi_u$ is
		\[
		\ell^2(\mathcal{G}_u)^{\pi_u}
		=
		\begin{cases}
			\mathbb C 1_{\mathcal{G}_u}, & |\mathcal{G}_u|<\infty,\\[2mm]
			\{0\}, & |\mathcal{G}_u|=\infty.
		\end{cases}
		\]
		Here $1_{\mathcal{G}_u}$ denotes the constant function of value one on $\mathcal{G}_u$.
	\end{lemma}

	\begin{lemma}\label{Le2} For each $u\in \mathcal{G}^{(0)}$, let
		\[
		U_u:\ell^2\bigl(\mathcal{H}_{\phi(u)}\bigr)\longrightarrow \ell^2(\mathcal{G}_u)
		\]
		be the unitary induced by the fibre bijection $\phi_u$. Then:
		\begin{enumerate}
			\item[\rm(i)] for all $f\in C_c(\mathcal{H})$,
			$
			U_u\pi_{\phi(u)}(f)=\pi_u(\phi^*f)U_u;
			$
			\item[\rm(ii)]
			$
			\Psi_\mathcal{G}(\phi^*f)=\phi^*(\Psi_\mathcal{H}(f));
			$
			\item[\rm(iii)]
			$
			U_u\ell^2\bigl(\mathcal{H}_{\phi(u)}\bigr)^{\pi_{\phi(u)}}
			=
			\ell^2(\mathcal{G}_u)^{\pi_u},
			$
			and hence $U_u$ also identifies the corresponding orthogonal complements.
		\end{enumerate}
	\end{lemma}

	\begin{proof}
		Part (i) follows from the functoriality of $\phi$ and the bijectivity of $\phi$ on source fibres. Indeed, for $\eta\in \mathcal{H}_{\phi(u)}$,
		\[
		U_u\pi_{\phi(u)}(f)\delta_\eta
		=
		\sum_{\xi\in \mathcal{H}_{\phi(u)}} f(\xi)\delta_{\phi_u^{-1}(\xi\eta)},
		\]
		while
		\[
		\pi_u(\phi^*f)U_u\delta_\eta
		=
		\sum_{\gamma\in \mathcal{G}_u} f(\phi(\gamma))\delta_{\gamma\phi_u^{-1}(\eta)}.
		\]
		Since $\phi$ is a functor and $\phi_u$ is a bijection, these two expressions agree. Equivalently, it is the \'etale counting-measure form of the regular-representation identity underlying \cite	[Proposition 2.7]{AustinMitra}.
		
		For (ii), using bijectivity on range fibres,
		\[
		\Psi_\mathcal{G}(\phi^*f)(u)
		=
		\sum_{\gamma\in \mathcal{G}^u} f(\phi(\gamma))
		=
		\sum_{\eta\in \mathcal{H}^{\phi(u)}} f(\eta)
		=
		\Psi_\mathcal{H}(f)(\phi(u)).
		\]
		Finally, (iii) follows from Lemma \ref{Le1}, the bijection $\phi_u$ preserves the cardinality of the source fibre and sends the constant function on $\mathcal{H}_{\phi(u)}$ to the constant function on $\mathcal{G}_u$.
	\end{proof}

	\begin{theorem}\label{Th1}
		Let $\mathcal{G}$ and $\mathcal{H}$ be \'etale groupoids with compact unit spaces, and let
		$
		\phi:\mathcal{G}\to \mathcal{H}
		$
		be a Haar system-preserving homomorphism. If 	$
		(\mathcal{H},\mathcal{R}_\mathcal{H})$ has Property {\rm(T)}, then
		$
		(\mathcal{G},\mathcal{R}_\mathcal{G})$ has Property {\rm(T)}.

	\end{theorem}
	
	\begin{proof}  Assume that $(K,c)$ is a Kazhdan pair for $	\mathcal{R}_\mathcal{H}$.  Since $\phi$ is proper and $K$ is compact, it follows that
		$
		L:=\phi^{-1}(K)
		$
		is compact. We next show that $
		(L,c)
		$
		is a Kazhdan pair for $	\mathcal{R}_\mathcal{G}$.

		Fix $u\in \mathcal{G}^{(0)}$ and $\xi\in \ell^2(\mathcal{G}_u)_{\pi_u}$. Set
		$
		\eta=U_u^*\xi.
		$
		By Lemma \ref{Le2}(iii),
		\[
		\eta\in \ell^2\bigl(\mathcal{H}_{\phi(u)}\bigr)_{\pi_{\phi(u)}}.
		\]
		Since $(K,c)$ is a Kazhdan pair for $\mathcal{R}_\mathcal{H}$, there exists $f\in C_c(\mathcal{H})$ such that
		\[
		\operatorname{supp}(f)\subseteq K,\qquad \|f\|_I\le 1,
		\]
		and
		\[
		\|\pi_{\phi(u)}(f)\eta-\pi_{\phi(u)}(\Psi_\mathcal{H}(f))\eta\|\ge c\|\eta\|.
		\]
		Put
		$
		\tilde f=\phi^*f.
		$
		Then $\operatorname{supp}(\tilde f)\subseteq \phi^{-1}(K)=L$, and $\|\tilde f\|_I\le \|f\|_I\le 1$. Moreover, by Lemma \ref{Le2},
		\[
		\pi_u(\tilde f)\xi-\pi_u(\Psi_\mathcal{G}(\tilde f))\xi
		=
		U_u\bigl(\pi_{\phi(u)}(f)\eta-\pi_{\phi(u)}(\Psi_\mathcal{H}(f))\eta\bigr).
		\]
		Taking norms gives
		\[
		\|\pi_u(\tilde f)\xi-\pi_u(\Psi_\mathcal{G}(\tilde f))\xi\|
		\ge c\|\eta\|=c\|\xi\|.
		\]
		Thus $(L,c)$ is a Kazhdan pair for $\mathcal{R}_\mathcal{G}$.
	\end{proof}

	\begin{remark}
		Surjectivity of $\phi$ is not required in Theorem \ref{Th1}. In \cite[Proposition 2.7]{AustinMitra} , surjectivity is used only to ensure that the induced map
		\[
		C^*_r(\mathcal{H})\longrightarrow C^*_r(\mathcal{G})
		\]
		is isometric. Consequently, regular Property (T) pulls back along   Haar system-preserving morphisms without any surjectivity hypothesis.
		
	\end{remark}

	\subsection{General representation families}
	
	The regular case works because the pullback morphism not only
	intertwines the regular representations but also identifies their
	constant-vector spaces. For an arbitrary representation $(H,\pi)$ of $C_c(\mathcal{G})$, every constant vector for $\pi$ is automatically constant for the restricted representation $\pi \circ \phi^{*}$, i.e., $H^\pi\subseteq H^{\pi\circ\phi^*}$. However, the reverse inclusion can fail, which is precisely the obstruction to	applying the Kazhdan inequality for the target family. Thus the passage from regular Property (T) to Property (T) for arbitrary representation 	families is not automatic. We have the following definition.
		
	\begin{definition}\label{De1}
Let $\phi:\mathcal{G}\to \mathcal{H}$	be a Haar system-preserving morphism. Let $\mathcal F_\mathcal{G}$ and $\mathcal F_\mathcal{H}$ be families of representations of $C_c(\mathcal{G})$ and $C_c(\mathcal{H})$, respectively. We say that $(\mathcal F_\mathcal{G},\mathcal F_\mathcal{H})$ is constant-vector compatible with $\phi$ if for every $(H,\pi)\in\mathcal F_\mathcal{G}$, the following conditions hold:
		\begin{enumerate}	\item[(i)] $\pi\circ\phi^*$ belongs to $\mathcal F_H$;
			\item[(ii)] the constant-vector spaces coincide, i.e.,
			$
			H^\pi=H^{\pi\circ\phi^*}.
			$
		\end{enumerate}
	\end{definition}
	Under this compatibility assumption, Kazhdan pairs pull back as expected.		
	\begin{theorem}\label{Th3}
		Assume $(\mathcal F_\mathcal{G},\mathcal F_\mathcal{H})$ is constant-vector compatible with $\phi$. If $(K,c)$ is a Kazhdan pair for $\mathcal F_\mathcal{H}$, then
		$
		(\phi^{-1}(K),c)
		$
		is a Kazhdan pair for $\mathcal F_\mathcal{G}$. Hence, if $(\mathcal{H},\mathcal F_\mathcal{H})$ has property {\rm(T)}, then $(\mathcal{G},\mathcal F_\mathcal{G})$ also has property {\rm(T)}.
	\end{theorem}
	
	\begin{proof}
		Fix $(H,\pi)\in\mathcal F_\mathcal{G}$ and $\xi\in H_\pi$. By Definition \ref{De1}(ii),
		$
		\xi\in H_{\pi\circ\phi^*}.
		$
		Choose $f\in C_c(\mathcal{H})$ from the Kazhdan property for $\mathcal F_\mathcal{H}$, and set
		$
		\tilde f=\phi^*f.
		$
		Since
		$
		\operatorname{supp}(\phi^*f)\subseteq \phi^{-1}(\operatorname{supp}(f))\subseteq \phi^{-1}(K)
		$
		and $\phi$ is proper, we have
		\[
		\operatorname{supp}(\tilde f)\subseteq \phi^{-1}(K).
		\] By the $I$-norm decreasing property of the pullback,
		$
		\|\tilde f\|_I\le \|f\|_I\le 1.
		$
		And by   the corresponding compatibility result for $\phi$, Haar system-preservation gives
		\[
		\Psi_\mathcal{G}(\phi^*f)= \phi ^*(\Psi_\mathcal{H}(f)).
		\]
		Consequently,
		$
		\|\pi(\tilde f)\xi-\pi(\Psi_\mathcal{G}(\tilde f))\xi\|
		=
		\|(\pi\circ\phi^*)(f)\xi-(\pi\circ\phi^*)(\Psi_\mathcal{H}(f))\xi\|
		\ge c\|\xi\|.
		$
	\end{proof}
	
	 	\begin{remark}
		For actions of a discrete group $\Gamma$ on a compact space $X$,
		Dell'Aiera--Willett obtain a more precise description: transformation groupoid $X\rtimes \Gamma$ has
		Property (T) precisely when $\Gamma$ has Property (T) with respect to
		the class of unitary representations arising in covariant representations
		of the action \cite[Proposition 4.19]{DellAiera Willett}. In particular, full Property (T) of
		$\Gamma$ implies full Property (T) of $X\rtimes \Gamma$. This conclusion
		relies on the specific representation theory of transformation groupoids
		and is stronger than what follows merely from the existence of a pullback
		$\ast$-homomorphism.
	\end{remark}

	\begin{remark}\label{cor3} Austin--Mitra define a partial morphism from $\mathcal{G}$ to $\mathcal{H}$ as a pair $(\phi,\mathcal{K})$, where $\mathcal{K}$ is an open   Haar subgroupoid   of $\mathcal{G}$ and $\phi: \mathcal{K}\to \mathcal{H}$ is a Haar system-preserving homomorphism \cite[Definition 3.3]{AustinMitra}. 
The preceding argument applies directly to the morphism $\phi: \mathcal{K}\to \mathcal{H}$ and yields that, if $\mathcal{G},\mathcal{H},\mathcal{K}$ have compact unit spaces and $(\mathcal{H},\mathcal{R}_\mathcal{H})$ has Property {\rm(T)}, then $(\mathcal{K},\mathcal{R}_\mathcal{K})$ has Property {\rm(T)}. Without additional hypotheses relating Property (T) of the open subgroupoid $\mathcal{K}$ to Property (T) of the ambient groupoid $\mathcal{G}$, no conclusion for $\mathcal{G}$ follows from this statement. \end{remark}

	\section{ Property (T) for transformation groupoids}

	Throughout this section, $\mathcal{G}$ is a locally compact Hausdorff \'etale groupoid with compact unit space $\mathcal{G}^{(0)}$, and $X$ is a compact Hausdorff left $\mathcal{G}$-space with (surjective) anchor map $p:X\to \mathcal{G}^{(0)}$. The associated transformation groupoid is
	\[
 	\mathcal{G}\ltimes X=\{(g,x)\in \mathcal{G}\times X: s(g)=p(x)\}.
	\]
	with
	$
	s(g,x)=x$, $r(g,x)=gx,
	$
	and   multiplication
	$
	(h,gx)(g,x)=(hg,x)
	$.
	The unit space $(\mathcal{G}\ltimes X)^{(0)}$ is identified with $X$ via
	$
	x\longmapsto (p(x),x),
	$
	and the inverse of $(g,x)\in\mathcal{G}\ltimes X$ is
	$
	(g,x)^{-1}=(g^{-1},gx).
	$
	Since $\mathcal{G}$ is \'etale, $\mathcal{G}\ltimes X$ is again a locally compact Hausdorff \'etale groupoid with compact unit space $X$.

	Let $q:\mathcal{G}\ltimes X\longrightarrow \mathcal{G}$ be the canonical projection defined by $q(\gamma,x)=\gamma$ for $(\gamma,x)\in \mathcal{G}\ltimes X$.	We study the passage of Property (T) from a groupoid $\mathcal{G}$ to its transformation groupoid $\mathcal{G}\ltimes X$, and show that Kazhdan pairs pull back along the projection.

	\begin{lemma}\label{Le3}
		The map $q$ is a Haar system-preserving groupoid homomorphism.
	\end{lemma}

	\begin{proof}
		It is immediate that $q$ is continuous and preserves the groupoid structure. Let $C\subseteq \mathcal{G}$ be compact. Then
		\[
		q^{-1}(C)=\{(g,x)\in C\times X:s(g)=p(x)\}.
		\]
		This set is closed in $C\times X$, since it is the preimage of the diagonal in $\mathcal{G}^{(0)}\times \mathcal{G}^{(0)}$ under the continuous map
		\[
		(g,x)\longmapsto (s(g),p(x)),
		\]
		and $\mathcal{G}^{(0)}$ is Hausdorff. As $C\times X$ is compact, $q^{-1}(C)$ is compact; hence $q$ is proper.
		
		Fix $x\in X$. Each $g\in \mathcal{G}_{p(x)}$ corresponds to a unique arrow $(g,x)\in(\mathcal{G}\ltimes X)_x$, so $q$ restricts to a bijection
		\[
		(\mathcal{G}\ltimes X)_x\longrightarrow \mathcal{G}_{p(x)}.
		\]
		The analogous assertion for range fibres follows by inversion. Since $\mathcal{G}$ is \'etale, fibrewise bijectivity is equivalent to preservation of the counting Haar systems.
	\end{proof}	
	
From Theorem \ref{Th1}, we have the following corollary.
	
	\begin{corollary}\label{Cor3}
		 If $(\mathcal{G},\mathcal{R}_\mathcal{G})$ has Property {\rm(T)}, then $(\mathcal{G}\ltimes X,\mathcal{R}_{\mathcal{G}\ltimes X})$ has Property {\rm(T)}.
		
		More precisely, if $(K,c)$ is a Kazhdan pair for $\mathcal{R}_\mathcal{G}$, then
		$
		(q^{-1}(K),c)
		$
		is a Kazhdan pair for $\mathcal{R}_{\mathcal{G}\ltimes X}$.
	\end{corollary}
	
		\begin{remark}  This corollary is a generalisation of \cite[Lemma 5.9]{BonickeDellAieraGabeWillett} and  no ampleness assumption is imposed. The proof of \cite[Lemma 5.9]{BonickeDellAieraGabeWillett} decomposes a Kazhdan function into finitely many pieces supported on compact open bisections, yielding a Kazhdan constant of the form $c/n$. The pullback argument above avoids this decomposition. Consequently, it removes the ample hypothesis and preserves the original constant $c$.
	\end{remark}

From Lemma \ref{le0}, for every representation $\pi: C_c(\mathcal{G})\to B(H)$, we have that $H^{\pi}\subseteq H^{\pi\circ q^{*}}$, where $q^*: C_c(\mathcal{G})\to C_c(\mathcal{G}\ltimes X)$ is the dual map of $q$ defined by $q^*(f)=fq$ for $f\in C_c(\mathcal{G})$. The next density statement is a well-known fact and is the technical point that allows one to show that $\pi$ and $\pi \circ q^{*}$ have the same constant vectors, i.e., $H^{\pi}=H^{\pi\circ q^*}$. For completion, we give its proof similar to   \cite[Proposition 1.72]{Goehle}.

	\begin{lemma}\label{Le4}
		The set
		\[
		A:=\operatorname{span}\{a q^*(f): a\in C(X),\ f\in C_c(\mathcal{G})\}
		\]
		is dense in $C_c(\mathcal{G}\ltimes X)$ with respect to the inductive-limit topology. Equivalently, for every $F\in C_c(\mathcal{G}\ltimes X)$, there are a compact set $D\subseteq\mathcal{G}\ltimes X$ and a net $\{F_i\}$ in $A$, all supported in $D$, such that
		$
		\|F_i-F\|_\infty\longrightarrow 0.
		$
		Consequently,
		$
		\|F_i-F\|_{I,\mathcal{G}\ltimes X}\longrightarrow 0.
		$
	\end{lemma}

	\begin{proof}
		Let $\mathcal{G}\, _r\ast_p X:=\{(g,y)\in\mathcal{G}\times X:r(g)=p(y)\}$ be the fibred product defined by $r$ and $p$.
Then the map $\Theta:(g,x)\in\mathcal{G}\ltimes X\longrightarrow (g,gx)\in \mathcal{G}_r\ast_p X$ is a homeomorphism with inverse $
		\Theta^{-1}(g,y)=(g,g^{-1}y)$ for $(g,y)\in \mathcal{G}\, _r\ast_p X$. Hence the map $$\Theta_*: f\in C_c(\mathcal{G}\ltimes X)\longrightarrow f\Theta^{-1}\in C_c(\mathcal{G}_r\ast_p X)$$ is a  $\ast$-isomorphism which is continuous uniformly.
		
 For $a\in C(X)$ and $f\in C_0(\mathcal{G})$, we let $f\otimes a\in C_0(\mathcal{G}\times X)$ defined by $(f\otimes a)(g,y)=f(g)a(y)$ for $(g,y)\in \mathcal{G}\times X$. Remark that the set $\operatorname{span}\{f\otimes a: a\in C(X), f\in C_0(\mathcal{G})\}$ is dense uniformly in $C_0(\mathcal{G}\times X)$. In particular, for $a\in C(X)$ and $f\in C_c(\mathcal{G})$, a direct computation of the convolution yields that $\Theta_*(aq^*(f))(g,y)=a(y)f(g)=(f\otimes a)(g,y)$ for $(g,y)\in \mathcal{G}_r\ast_p X$.

Given $F\in C_c(\mathcal{G}\ltimes X)$, we have $\Theta_*(F)\in C_c(\mathcal{G}\, _r\ast_p X)$ and extend it a function in $C_0(\mathcal{G}\times X)$, also denoted by $\Theta_*(F)$.  Therefore it follows from the above paragraph that we can find a net $\{k_i\}$ converges to $\Theta_*(F)$ uniformly, where $k_i=\sum\limits_{j=1}^{m_i}f_{ij}\otimes a_{ij}$ for $\{f_{ij}\}\subset C_0(\mathcal{G})$ and $\{a_{ij}\}\subset C(X)$.

Let $K=\operatorname{supp}(\Theta_*(F))$ and $L=q(\operatorname{supp}(F))$. Then $L$ is compact in $\mathcal{G}$. 	Choose $\chi\in C_c(\mathcal{G})$ such that $0\leq \chi\leq 1$ and $\chi(g)=1$ for each $g\in L$. Put $C=\operatorname{supp}(\chi)$ and $D=\{(g,y)\in \mathcal{G}\, _r\ast_p X:\, g\in C\}$. Then both of $C$ and $D$ are compact, $L\subseteq C$ and $K\subseteq D$. For each $i$, let $\widetilde{k_i}=\sum\limits_{j=1}^{m_i}(f_{ij}\cdot \chi)\otimes a_{ij}$, where $f_{ij}\cdot \chi$ denotes the point-wise product in $C_0(\mathcal{G})$. We can check that the support of each $\widetilde{k_i}$ is contained $D$ and $\{\widetilde{k_i}\}$ converges to $\Theta_*(F)$ uniformly. Hence each $\Theta_*^{-1}(\widetilde{k_i})$ has support contained in $\Theta^{-1}(D)$ and $\{\Theta_*^{-1}(\widetilde{k_i})\}$ converges to $F$ uniformly.\end{proof}
		

	Using the density result, we can now compare the constant vectors for $\pi$ and $\pi\circ q^*$.
	\begin{proposition}\label{Pro1}
		Let
		$
		\pi:C_c(\mathcal{G}\ltimes X)\longrightarrow B(  H)
		$
		be a representation.
		Then
		$H^\pi=  H^{\pi\circ q^*}$, thus $H_\pi=  H_{\pi\circ q^*}$.
	\end{proposition}
	
	\begin{proof} Put
		$
		\rho:=\pi\circ q^*:C_c(\mathcal{G})\longrightarrow B(  H).
		$
		From Lemma \ref{le0}, we have $H^{\pi}\subseteq H^\rho$. Conversely, let $\xi\in  H^\rho$. If $a\in C(X)$ and $f\in C_c(\mathcal{G})$, then
		\[
		\pi(a*q^*(f))\xi
		=\pi(a)\pi(q^*(f))\xi
		=\pi(a)\pi((\Psi_\mathcal{G}(f))\circ p)\xi.
		\]
		On the other hand, for $y\in X$,
		\[
		\begin{aligned}
			\Psi_{\mathcal{G}\ltimes X}(a*q^*(f))(y)
			=\sum_{\gamma\in (\mathcal{G}\ltimes X)_y}(a*q^*(f))(\gamma)
			=a(y)\sum_{g\in \mathcal{G}_{p(y)}}f(g)
			=a(y)\Psi_\mathcal{G}(f)(p(y)).
		\end{aligned}
		\]	Therefore,
		$$
		\pi(a*q^*(f))\xi=\pi(\Psi_{\mathcal{G}\ltimes X}(a*q^*(f)))\xi.
		$$
		By linearity the same holds for every element of $A$ from Lemma \ref{Le4}. Since representations are $I$-norm bounded and
		\[
		\|\Psi_{\mathcal{G}\ltimes X}(F)\|_\infty\le \|F\|_{I,\mathcal{G}\ltimes X},
		\]
		Lemma \ref{Le4} allows us to pass to the limit and conclude that
		\[
		\pi(F)\xi=\pi(\Psi_{\mathcal{G}\ltimes X}(F))\xi
		\]
		for all $F\in C_c(\mathcal{G}\ltimes X)$. Hence $\xi\in  H^\pi$. Therefore $  H^\pi=  H^\rho$.
	\end{proof}

	We now formulate the main pullback theorem for general representation families.
	
	\begin{theorem}\label{Th2}
		Let $\mathcal F$ be a family of representations of $C_c(\mathcal{G})$ and let $\mathcal E$ be a family of representations of $C_c(\mathcal{G}\ltimes X)$. Assume that for every $\pi\in\mathcal E$, the restricted representation
		$
		\pi\circ q^*
		$
		belongs to $\mathcal F$  up to unitary equivalence. If $(\mathcal{G},\mathcal F)$ has property {\rm(T)}, then $(\mathcal{G}\ltimes X,\mathcal E)$ has property {\rm(T)}.
		
		More precisely, if $(K,c)$ is a Kazhdan pair for $(\mathcal{G},\mathcal F)$, then $(L,c)$ is a Kazhdan pair for $(\mathcal{G}\ltimes X,\mathcal E)$, where $L:=q^{-1}(K)=\{(g,x)\in \mathcal{G}\ltimes X:g\in K\}$.
		\end{theorem}
	
	\begin{proof} Since $q$ is proper, the set $L$ is compact in $\mathcal{G}\ltimes X$. 		
		Fix $\pi\in\mathcal E$ and
		$
		\xi\in   H_\pi .
		$
		Put
		$
		\rho:=\pi\circ q^*.
		$
		By Proposition \ref{Pro1},
		$
		\xi\in  H_\rho.
		$
		Since $(K,c)$ is a Kazhdan pair for $(\mathcal{G},\mathcal F)$, there exists $f\in C_c(\mathcal{G})$ such that
		\[
		\operatorname{supp}(f)\subseteq K,
		\qquad
		\|f\|_{I,\mathcal{G}}\le 1,
		\]
		and
		\[
		\|\rho(f)\xi-\rho(\Psi_\mathcal{G}(f))\xi\|\ge c\|\xi\|.
		\]
		Set
		$
		F:=q^*(f).
		$
		Then
		\[
		\operatorname{supp}(F)\subseteq L,
		\qquad
		\|F\|_{I,\mathcal{G}\ltimes X}\le 1,
		\]
		and
		\[
		\Psi_{\mathcal{G}\ltimes X}(F)=(\Psi_\mathcal{G}(f))\circ p.
		\]
		Therefore
		\[
		\begin{aligned}
			\|\pi(F)\xi-\pi(\Psi_{\mathcal{G}\ltimes X}(F))\xi\|
			=
			\|\rho(f)\xi-\rho(\Psi_\mathcal{G}(f))\xi\|
			\ge c\|\xi\|.
		\end{aligned}
		\]
		Thus $(L,c)$ is a Kazhdan pair for $(\mathcal{G}\ltimes X,\mathcal E)$.
	\end{proof}

		Taking all representations on both sides gives the following consequence.
	\begin{corollary}\label{Cor1}
		If $\mathcal{G}$ has topological Property {\rm(T)}  with respect to all representations of $C_c(\mathcal{G})$, then $\mathcal{G}\ltimes X$ has topological Property {\rm(T)}.
	\end{corollary}

	The same argument applies to representations factoring through the reduced algebra.
	\begin{corollary}
		Let $\operatorname{Rep}_r(\mathcal{G})$ denote the family of representations of $C_c(\mathcal{G})$ that factor through $C^*_r(\mathcal{G})$, and define $\operatorname{Rep}_r(\mathcal{G}\ltimes X)$ similarly. If $(\mathcal{G},\operatorname{Rep}_r(\mathcal{G}))$ has property {\rm(T)}, then $(\mathcal{G}\ltimes X,\operatorname{Rep}_r(\mathcal{G}\ltimes X))$ has property {\rm(T)}.
	\end{corollary}

	\begin{proof}
		For $f\in C_c(\mathcal{G})$,
		\[
		\begin{aligned}
			\|q^*(f)\|_{C^*_r(\mathcal{G}\ltimes X)}
			&=\sup_{x\in X}\|\pi_x(q^*(f))\|
			=\sup_{x\in X}\|\pi_{p(x)}(f)\|\\
			&\le\sup_{u\in \mathcal{G}^{(0)}}\|\pi_u(f)\|
			=\|f\|_{C^*_r(\mathcal{G})}.
		\end{aligned}
		\]
		Hence $q^*$ extends to a $*$-homomorphism
		\[
		q^*_r:C^*_r(\mathcal{G})\longrightarrow C^*_r(\mathcal{G}\ltimes X).
		\]
		Now let $\pi\in\operatorname{Rep}_r(\mathcal{G}\ltimes X)$. By definition, $\pi$ factors through $C^*_r(\mathcal{G}\ltimes X)$, so there exists a representation $\pi_r$ of $C^*_r(\mathcal{G}\ltimes X)$ such that $\pi=\pi_r\circ q_{\mathcal{G}\ltimes X}$, where $q_{\mathcal{G}\ltimes X}:C_c(\mathcal{G}\ltimes X)\to C^*_r(\mathcal{G}\ltimes X)$ is the natural map. Since $q_{\mathcal{G}\ltimes X}\circ q^*=q^*_r\circ q_\mathcal{G}$, we obtain
		\[
		\pi\circ q^*=\pi_r\circ q^*_r\circ q_\mathcal{G},
		\]
		so $\pi\circ q^*$ factors through $C^*_r(\mathcal{G})$, i.e.\ $\pi\circ q^*\in\operatorname{Rep}_r(\mathcal{G})$. Therefore Theorem \ref{Th2} applies with $\mathcal F=\operatorname{Rep}_r(\mathcal{G})$ and $\mathcal E=\operatorname{Rep}_r(\mathcal{G}\ltimes X)$.
	\end{proof}
	
	The compatibility hypothesis in Theorem~\ref{Th2} cannot be omitted in general, as the following example illustrates.

	\begin{example}\label{Ex1}
		
		Let $\mathcal{G}=\mathbb Z$, viewed as a groupoid with one unit, and let $X=\{*\}$ with the trivial action. Then
		$$
		\mathcal{G}\ltimes X\cong \mathbb Z.
		$$
		Let
		$
		\mathcal F=\{1\},
		$
		where $1$ is the trivial representation of $\mathbb Z$. Since
		$
	  H_1 =\{0\},
		$
		$\mathbb Z$ has property (T) with respect to $\mathcal F$ vacuously.

		On the other hand, let
		$
		\mathcal E=\{\lambda\},
		$
		where $\lambda$ is the left regular representation on $\ell^2(\mathbb Z)$. We claim that $\mathbb Z$ does not have property (T) with respect to $\mathcal E$.

		Indeed, since $\mathbb Z$ is infinite,
		$
		  H^\lambda=\{0\},
		$
		so
		$
		  H_\lambda =\ell^2(\mathbb Z).
		$
		For $N\ge 1$, define
		\[
		\xi_N:=\frac{1}{\sqrt{2N+1}}\mathbf 1_{[-N,N]}.
		\]
		If $K\subseteq \mathbb Z$ is finite, then
		\[
		\max_{k\in K}\|\lambda_k\xi_N-\xi_N\|\longrightarrow 0.
		\]
		For any $f\in C_c(\mathbb Z)$ with
		$
		\operatorname{supp}(f)\subseteq K$,
		$
		\|f\|_1\le 1,
		$
		we have
		\[
		\begin{aligned}
			\|\lambda(f)\xi_N-\lambda(\Psi_{\mathbb Z}(f))\xi_N\|
			&=\left\|\sum_{k\in K}f(k)(\lambda_k\xi_N-\xi_N)\right\|\\
			&\le \max_{k\in K}\|\lambda_k\xi_N-\xi_N\|\sum_{k\in K}|f(k)|\\
			&\le \max_{k\in K}\|\lambda_k\xi_N-\xi_N\|.
		\end{aligned}
		\]
		The right-hand side tends to zero as $N\to\infty$. Hence no constant $c>0$ can satisfy the Kazhdan inequality for all $N$, so $\mathbb Z$ does not have property (T) with respect to $\mathcal E$. This shows that the compatibility hypothesis in Theorem \ref{Th2} is necessary.

	\end{example}


	The direction in Theorem \ref{Th1} is contravariant: Property (T) on the target implies Property (T) on the source. The converse fails, even under surjectivity and even for full Property (T).
	
	\begin{proposition}
		There exists a surjective   continuous Haar system-preserving morphism
		$
		\phi:\mathcal{G}\longrightarrow \mathcal{H}
		$
		between locally compact Hausdorff \'etale groupoids with compact unit spaces such that $\mathcal{G}$ has full topological Property {\rm(T)} but $\mathcal{H}$ does not.
	\end{proposition}
	
	\begin{proof}
		Let $G=F_2$ and let $X=\partial F_2$ be its Gromov boundary. Set
		\[
		\mathcal{G}=X\rtimes F_2,\qquad \mathcal{H}=F_2,
		\]
		where $F_2$ is regarded as a one-unit groupoid. Dell'Aiera--Willett \cite[Example 6.11]{DellAiera Willett} show that $X\rtimes F_2$ has topological Property (T). The free group $F_2$ does not have Kazhdan's Property (T).
		
		The projection
		\[
		q:X\rtimes F_2\longrightarrow F_2,\qquad q(x,\gamma)=\gamma,
		\]
		is a surjective continuous groupoid morphism. Since $X$ is compact, $q$ is proper. On each source fibre it is a bijection, and likewise on each range fibre; hence it preserves the counting Haar systems. Thus $q$ satisfies all the morphism properties in the statement, but Property (T) does not pass from $\mathcal{G}$ to $\mathcal{H}$.
	\end{proof}
	
	\begin{remark}
		By \cite[Proposition 2.7]{AustinMitra}, surjectivity of $q$ implies that the pullback map
		\[
		q^*:C^*_r(\mathcal{H})\longrightarrow C^*_r(\mathcal{G})
		\]
		is isometric. Therefore even an isometric pullback embedding of reduced groupoid $C^*$-algebras does not imply that Property (T) descends from the source groupoid to the target groupoid. The obstruction is representation-theoretic rather than merely norm-theoretic.
	\end{remark}

\section{Ampliation invariance of Property (T)}

	There are some notions of equivalence of groupoids, such as Morita equivalence, groupoid equivalence, equivalence via a linking
groupoid, equivalence via isomorphic ampliations, Kakutani equivalence. From \cite[Proposition 3.10, Theorem 3.12]{Farsi}, the first four notions coincide for \'{e}tale groupoids, and they all coincide for ample \'etale groupoids with $\sigma$-compact unit spaces.

Given an \'{e}tale groupoid $\mathcal{G}$ and a surjective local homeomorphism $\varphi$ from another locally compact Hausdorff space $X$ onto $\mathcal{G}^{(0)}$, we define the ampliation $\mathcal{G}^{\varphi}$ of $\mathcal{G}$ with respect to $\varphi$ as
	\[
	\mathcal{G}^{\varphi}:=\{(x,\gamma,y)\in X\times 	\mathcal{G} \times X:\varphi(x)=r(\gamma),\ \varphi(y)=s(\gamma)\}.
	\]
One can check that it is an \'{e}tale groupoid with multiplication $
	(x,\gamma,y)(y,\eta,z)=(x,\gamma\eta,z)$ for $(\gamma,\eta)\in \mathcal{G}^{(2)}$ and inverse map $(x,\gamma,y)^{-1}=(y,\gamma^{-1}, x)$ when  endowed with the relative product topology from $X\times \mathcal{G}\times X$. By identifying $(x,\varphi(x), x)$ with $x$, we can identify $(\mathcal{G}^{\varphi})^{(0)}$ with $X$.

 For our purpose, we say that two \'{e}tale groupoids $\mathcal{G}$ and $\mathcal{H}$ are Morita equivalent if they admits isomorphic ampliations, and say that $\mathcal{G}$ and $\mathcal{H}$ are Kakutani equivalent if there are full clopen subsets $X\subseteq \mathcal{G}^{(0)}$ and $Y\subseteq \mathcal{H}^{(0)}$ such that the restriction $\mathcal{G}|_X$ and $\mathcal{H}|_Y$ are isomorphic.

	In this  section, we show that each ample groupoid (with compact unit space) and its ampliation groupoid (with compact unit space) have the same property (T). We assume that $\mathcal{G}$ is an ample \'{e}tale groupoid with compact unit space and $\varphi$ is a surjective local homeomorphism from a compact Hausdorff space $X$ onto $\mathcal{G}^{(0)}$. We first record the following useful fact on local homeomorphisms and then, use it to obtain a matrix-corner description of the associated function algebra  $C_c(\mathcal{G}^{\varphi})$.

	\begin{lemma}\label{Le6}
		Let $X$ and $Y$ be two compact Hausdorff spaces, let $Y$ be zero-dimensional and let
		$
		\varphi :X\to Y
		$
		be a surjective local homeomorphism.  Then the map, defined by $d(y):=\#\varphi^{-1}(y)$ for $y\in Y$, is continuous on $Y$.

Moreover, there exist an integer $N\ge 1$, pairwise disjoint clopen subsets
		$
		X_1,\dots,X_N\subseteq X,
		$
		and clopen subsets
		$
		U_1,\dots,U_N\subseteq Y
		$
		such that
		$$ U_1=Y,\;\;
		X=\bigsqcup_{i=1}^N X_i,
		\;\; \varphi|_{X_i}:X_i\to U_i$$ is a homeomorphism for every $i$, and
	$d(y)=\sum_{i=1}^N \mathbf 1_{U_i}(y)$.

	\end{lemma}

	\begin{proof}
		From the assumption, we have that $\varphi$ is a finite-sheeted covering map. Hence, for each $y\in Y$, $\varphi^{-1}(y)$ is finite and there exist an open subset $U_y\subseteq Y$ containing $y$ and disjoint open subsets $V_1,V_2,\cdots, V_n$ of $X$ such that $\varphi|_{V_i}: V_i\rightarrow U_y$ is a homeomorphism for each $i$ and $\varphi^{-1}(U_y)=\cup_{i=1}^nV_i$, where $n=d(y)$. Since $Y$ is zero-dimensional, $X$ has a base consisting of clopen subsets of $X$. Hence we may assume that $U_y$ and $V_i's$ are clopen. Since $Y$ is compact and zero-dimensional, there is a finite clopen partition
		\[
		Y=C_1\sqcup\cdots\sqcup C_m
		\]
		such that each $\varphi^{-1}(C_k)$ can be written as a disjoint union of clopen subsets of $X$:
		\[
		\varphi^{-1}(C_k)=C_{k,1}\sqcup\cdots\sqcup C_{k,d_k},
		\]
		where $\varphi|_{C_{k,i}}: C_{k,i}\rightarrow C_k$ is a homeomorphism for each $k,i$.
		
		Let $N=\max_k d_k$. For each $1\leq i\leq N$, we define
		\[
		X_i:=\bigsqcup_{\{k:\,d_k\ge i\}} C_{k,i},
		\qquad
		U_i:=\bigcup_{\{k:\,d_k\ge i\}} C_k.
		\]
		Then $X_i$ and $U_i$ are clopen, $\varphi|_{X_i}:X_i\to U_i$ is a homeomorphism, and $U_1=Y$. The formula for $d$ follows immediately.
	\end{proof}
	
	In the rest, we let $\mathcal{G}$, $X$, $\varphi$ and $\mathcal{G}^{\varphi}$ as before, and adopt the notations for $Y=\mathcal{G}^{(0)}$. Let $\sigma_i$ be the inverse of the restriction $\varphi|_{X_i}$ for each $i$. Then it is a homeomorphism from $U_i$ onto $X_i$. Hence we can rewrite the ampliation groupoid as
	 $$\mathcal{G}^{\varphi}=\{(\sigma_i(r(\gamma)), \gamma,\sigma_j(s(\gamma)):\, \gamma\in 	\mathcal{G}^{U_i}_{U_j}, 1\leq i,j\leq N\}.$$

 For the clopen cover $\mathcal{U}:=\{U_1,U_2,\cdots,U_N\}$ of $\mathcal{G}^{(0)}$, we recall that the localization of $\mathcal{G}$ over $\mathcal{U}$ is the \'{e}tale groupoid $$\mathcal{G}_{\mathcal{U}}=\{(i,\gamma,j): \gamma\in \mathcal{G}^{U_i}_{U_j}, 1\leq i,j\leq N\},$$
which is endowed with the relative product topology on $\{1,2,\cdots,N\}\times \mathcal{G}\times \{1,2,\cdots,N\}$, with multiplication $
	(i,\gamma,j)(j,\eta,k)=(i,\gamma\eta,k)$ for $(\gamma,\eta)\in \mathcal{G}^{(2)}$, and inverse $(i,\gamma,j)^{-1}=(j,\gamma^{-1}, i)$. Clearly, $\mathcal{G}_{\mathcal{U}}^{(0)}=\{(i,u,i):\, u\in U_i\}$.

The map $(i,\gamma,j)\rightarrow (\sigma_i(r(\gamma)),\gamma,\sigma_j(s(\gamma)))$ is a $*$-isomorphism preserving $\|\cdot\|_I$-norm from $\mathcal{G}_{\mathcal{U}}$ onto $\mathcal{G}^{\varphi}$. In the isomorphism, we identity $\mathcal{G}_{\mathcal{U}}^{(0)}$ with $X$ and let $\mathcal{G}^{\varphi}=\mathcal{G}_{\mathcal{U}}$.

Let $M_n(C_c(\mathcal{G}))$ be the algebra of all $n\times n$-matrices with entries in $C_c(\mathcal{G})$ and let
	\[
	P:=\operatorname{diag}(1_{U_1},\dots,1_{U_N})
	\]
	be the diagonal matrix in $M_n(C_c(\mathcal{G}))$. The following gives a matrix representation of $C_c(\mathcal{G}^{\varphi})$.

	\begin{proposition}\label{Pro2}
		There is a canonical $*$-isomorphism
		$		\Phi:C_c(\mathcal{G}^{\varphi})\to P M_N(C_c(\mathcal{G}))P $,
		defined by $\Phi(f)=(f_{ij})_{i,j=1}^N$ for $f\in C_c(\mathcal{G}^{\varphi})$, where
		$$f_{ij}(\gamma)=\left\{\begin{array}{ll}
		f(i,\gamma,j), & \mbox{ for $\gamma\in \mathcal{G}^{U_i}_{U_j}$}\\ 0,  & \mbox{ for otherwise,} \end{array}\right.$$
		for each $i,j$. Under this identification,
		\begin{equation}\label{eq1}
			\Psi_{\mathcal{G}^{\varphi}}(f)(i,u,i)=\sum_{j=1}^N \Psi_\mathcal{G}(f_{ij})(u),\qquad \mbox{ for $u\in U_i$},
		\end{equation}

		\[ \Phi(\Psi_{\mathcal{G}^{\varphi}}(f))=diag(\sum_{j=1}^N \Psi_\mathcal{G}(f_{1j}), \sum_{j=1}^N \Psi_\mathcal{G}(f_{2j}), \cdots, \sum_{j=1}^N \Psi_\mathcal{G}(f_{nj}))\]
and
		\[
		\|F\|_{I,\mathcal{G}^\varphi}
		=
		\max\left\{
		\sup_{1\le i\le N}\sup_{u\in U_i}
		\sum_{j=1}^N\sum_{h\in \mathcal{G}^u}|f_{ij}(h)|,\,
		\sup_{1\le j\le N}\sup_{u\in U_j}
		\sum_{i=1}^N\sum_{h\in \mathcal{G}_u}|f_{ij}(h)|
		\right\}.
		\]
	\end{proposition}

	\begin{proof} We can check that $\Phi$ is a $\ast$-isomorphism with inverse $\Phi^{-1}$, defined by
$$\Phi^{-1}((f_{ij}))(i,\gamma,j)=f_{ij}(\gamma)$$ for $(f_{ij})\in P M_N(C_c(\mathcal{G}))P$ and $(i,\gamma,j)\in \mathcal{G}^{\varphi}$.

		Every arrow of $\mathcal{G}^\varphi$ is of the form
		\[
		(\sigma_i(r(h)),h,\sigma_j(s(h)))
		\]
		for uniquely determined indices $i,j$ satisfying $r(h)\in U_i$ and $s(h)\in U_j$. This gives the asserted linear bijection. The formulas for multiplication and involution in $\mathcal{G}^\varphi$ translate exactly into matrix multiplication with convolution in the entries, proving that $\Phi$ is a $*$-isomorphism.

If let $f=\Phi^{-1}((f_{ij}))$, then $f(i,\gamma,j)=f_{ij}(\gamma)$ for $(i,\gamma,j)\in \mathcal{G}^{\varphi}$.
		
		Equation~\eqref{eq1} follows by summing over the range fibre of $\sigma_i(u)$. The equation for the $I$-norm follows similarly from the explicit descriptions of the range and source fibres.
	\end{proof}


The matrix-corner description of Proposition~\ref{Pro2} allows us to pass from representations of $C_c(\mathcal{G})$ to representations of $C_c(\mathcal{G}^\varphi)$ and back. For a representation $\rho:C_c(\mathcal{G})\to B(H)$, we have the representation $\pi_{\rho}: M_n(C_c(\mathcal{G}))\to B(\bigoplus_{i=1}^N   H)$ by $$\pi_{\rho}((f_{ij}))=(\rho(f_{ij})) \hbox{ for $(f_{ij})\in M_n(C_c(\mathcal{G}))$}.$$
 We write vectors $\oplus_{i=1}^N\xi_i$ in $\bigoplus_{i=1}^N   H$ as column vectors.  Let $K_{\rho}=\pi_{\rho}(P)(\bigoplus_{i=1}^N   H)$, i.e.,
$K_{\rho}=\bigoplus_{i=1}^N   H_i$ for $  H_i:=\rho(1_{U_i}) (H)$ ($i=1,2,\cdots,n$). Clearly, $\pi_{\rho}$ induces a representation of $PM_n(C_c(\mathcal{G}))P$ on $K_{\rho}$. We still denote the induced representation  by $\pi_{\rho}$. By the identification in Proposition \ref{Pro2}, $\pi_{\rho}$ is also a representation of $C_c(\mathcal{G}^\varphi)$ on $K_{\rho}$.

	\begin{lemma}\label{Le7}
		Every nondegenerate representation of $C_c(\mathcal{G}^\varphi)$ is unitarily equivalent to $\pi_\rho$ for a nondegenerate representation $\rho$ of $C_c(\mathcal{G})$. The representation $\rho$ is unique up to unitary equivalence.
	\end{lemma}
	
	\begin{proof} We let $\mathcal{B}=P\hbox{M}_N(C_c(\mathcal{G}))P$. Let $\{E_{ij}: 1\leq i,j\leq N\}$ be the canonical matrix unit system of $\hbox{M}_N(\mathbb{C})$, and for $a\in C_c(\mathcal{G})$, let $aE_{ij}$ be the matrix in $\hbox{M}_N(C_c(\mathcal{G}))$ whose entry is $a$ at $(i,j)$ and $0$ at other position. Hence each matrix $(a_{ij})$ in $\hbox{M}_N(C_c(\mathcal{G}))$ can be written as $\sum_{i,j=1}^N a_{ij}E_{ij}$. We note that $1_{U_1}$ is the unit in $C_c(\mathcal{G})$ and $P=\sum_{i=1}^N 1_{U_i}E_{ii}$ is the unit in $\mathcal{B}$.

Let $\pi$ be a unital representation of $\mathcal{B}$ on $K$. Put
		$H:=\pi(1_{U_1}E_{11})(K)$. Note that $\pi(aE_{11})(H)\subseteq H$ for each $a\in C_c(\mathcal{G})$. Then we have a unital representation $\rho: C_c(\mathcal{G})\rightarrow B(H)$, defined by  $\rho(a)=\pi(aE_{11})|_H$. As before, let
$K_{\rho}=\bigoplus_{i=1}^N H_i$ for $H_i:=\rho(1_{U_i})(H)=\pi(1_{U_i}E_{11})(H)$ ($i=1,2,\cdots,N$) and $\pi_{\rho}$ is the associated representation of $\mathcal{B}$ on $K_{\rho}$. If let $v_i=1_{U_i}E_{i1}$, then it is a partial isometry such that $v_iv_i^*=1_{U_i}E_{ii}$ and $v_i^*v_i=1_{U_i}E_{11}$ for each $i$. Define $$U: K_{\rho}\rightarrow K, \,\,\,\,\,\,\,\, U(\oplus_{i=1}^N\xi_i)=\sum_{i=1}^N\pi(v_i)\xi_i.$$
Then $U$ is unitary and intertwines $\pi_\rho$ with $\pi$.
	\end{proof}
	
	Next we fix a (unital) representation $\rho:C_c(\mathcal{G})\to B(  H)$. Define
	$$S_\rho:  H\to \mathcal L, \,\,\, S_\rho\xi:=\oplus_{i=1}^N\rho(1_{U_i})\xi \mbox{ for $\xi\in H$}.$$
Then $S_\rho^*S_\rho=\rho(d)$, where $d=\sum_{i=1}^N 1_{U_i}\in C(\mathcal{G}^{(0)})$. Let $S_{\rho}=W_{\rho}|S|=W_{\rho}\rho(d^{\frac 12})$ be the polar decomposition of $S_{\rho}$. Note that $d$ is inverse in $C_c(\mathcal{G}^{(0)})$. So $W_{\rho}: H\to \mathcal{L}$ is an isometry with range $S_{\rho}(H)$.
The following result characterizes the constant vectors of $\pi_\rho$.

	\begin{proposition}\label{Pro7}
		 $\mathcal L^{\pi_\rho}=S_\rho(  H^\rho).$ 

	\end{proposition}
	
	\begin{proof}
		Suppose first that 	$\xi=\oplus_{i=1}^N\xi_i\in \mathcal L^{\pi_\rho}.$
		For each $i$, consider the element $F$ of the matrix corner having $1_{U_i}$ in the $(i,1)$ entry and zero elsewhere. Its constant-vector equation, $\pi_{\rho}(F)\xi=\pi_{\rho}(\Psi_{\mathcal{G}^{\varphi}}(F))\xi$, gives
		$	\xi_i=\rho(1_{U_i})\xi_1.
		$
		Since $U_1=Y$, the $(1,1)$ corner contains all of $C_c(\mathcal{G})$. Applying the constant-vector equation to an element having arbitrary $f\in C_c(\mathcal{G})$ in the $(1,1)$ entry and zero elsewhere gives $
		\rho(f)\xi_1=\rho(\Psi_\mathcal{G}(f))\xi_1$.
		Thus $\xi_1\in  H^\rho$, and $
		\xi=S_\rho(\xi_1)\in S_{\rho}(H^{\rho}).$
		
		Conversely, let $\eta\in  H^\rho$ and put $\xi=S_\rho\eta$. For $F=(f_{ij})\in PM_n(C_c(\mathcal{G}))P$, since  $\sum\limits_{j=1}^N\rho(f_{ij})\eta=\sum\limits_j \rho(\Psi_{\mathcal{G}}(f_{ij}))\eta$ and $f_{ij}1_{U_j}=f_{ij}$ for each $i,j$, it follows from Proposition \ref{Pro2} that $\pi_{\rho}(F)\xi=\pi_{\rho}(\Psi_{\mathcal{G}^{\varphi}}(F))\xi$, thus $\xi\in \mathcal L^{\pi_\rho}$.\end{proof}

Let $q: \mathcal{G}^{\varphi} \rightarrow \mathcal{G}$ be the projection defined by $$q(i,\gamma,j)=\gamma \hbox{ for $(i,\gamma,j)\in \mathcal{G}^{\varphi}$}.$$ Then $q$ is a  proper surjective continuous homomorphism. 
Hence, for a compact subset $K\subseteq \mathcal{G}$,  the set $\widetilde{K}:=q^{-1}(K)$ is compact in $\mathcal{G}^{\varphi}$. In particular, $\widetilde{\mathcal{G}^{(0)}}=\{(i,u,j): u\in U_i\cap U_j\}$ is open and compact in $\mathcal{G}^{\varphi}$. 
Let $q^*: C_c(\mathcal{G})\rightarrow C_c(\mathcal{G}^{\varphi})$ be  the dual map, defined by $q^*(f)=fq$ for $f\in C_c(\mathcal{G})$. Recall that $\varphi^*(f)=f\varphi\in C(X)$ for each $f\in C(\mathcal{G}^{(0)})$ and 
that $d(y)$ is inverse in $C(\mathcal{G}^{(0)})$. One can check that $q^*(1_{\mathcal{G}^{0}})=1_{\widetilde{\mathcal{G}^{(0)}}}$ and
 $$e_{\varphi}:=\varphi^*{(d^{-1})}1_{\widetilde{\mathcal{G}^{(0)}}}$$
 is a projection in $C_c(\mathcal{G}^{\varphi})$ and $\Psi_{\mathcal{G}^{\varphi}}(e_{\varphi})(i,u,i)=1$ for all $u\in U_i$, i.e.,  $\Psi_{\mathcal{G}^{\varphi}}(e_{\varphi})$ is the identity element $1_X$ in $C_c(\mathcal{G}^{\varphi})$. We can also claim that $e_{\varphi}$ has the matrix representation: $$\Phi(e_{\varphi})=(d^{-1}1_{U_i}1_{U_j})_{ij}$$ and $\|e_\varphi\|_{I,\mathcal{G}^\varphi}=1$. Thus $\pi_{\rho}(e_{\varphi})=(\rho(d^{-1})\rho(1_{U_i})\rho(1_{U_j}))_{ij}$. A direct calculation shows that $\pi_{\rho}(e_{\varphi})$ has the range $S_{\rho}(H)$, thus it is the range projection of $W_{\rho}$, i.e., $$\pi_\rho(e_\varphi)=W_\rho W_\rho^*.$$


 \begin{lemma}\label{Le9}  For $f\in C_c(\mathcal{G})$, define $J(f)\in C_c(\mathcal{G}^{\varphi})$ by
	$$J(f)=\varphi^*(d^{-\frac 12})q^*(f)\varphi^*(d^{-\frac 12}).$$
Then $\operatorname{supp}(J(f))\subseteq \widetilde{\operatorname{supp}(f)}$, $\Psi_{\mathcal{G}^{\varphi}}(J(f))=\varphi^*(\Psi_{\mathcal{G}}(d^{-\frac 12} f d^{\frac 12} ))$ and $\|J(f)\|_{I,\mathcal{G}^\varphi}\le \sqrt{N}\|f\|_{I,\mathcal{G}}$ for each $f$. Moreover,
		\begin{equation}
			\pi_\rho(J(f))=W_\rho\rho(f)W_\rho^*, \mbox{ and } W_\rho^*\pi_\rho(\Psi_{\mathcal{G}^{\varphi}}(J(f)))W_\rho=\rho(\Psi_{\mathcal{G}}(fd^{\frac 12}))\rho(d^{-\frac 12}).
			\label{eq4}
		\end{equation}
	\end{lemma}
	
	\begin{proof} The support statement is immediate. For $f\in C_c(\mathcal{G})$, the element $J(f)\in C_c(\mathcal{G}^\varphi)$ has the matrix representation:
	$$J(f)_{ij}:=(1_{U_i}d^{-\frac 12})f(1_{U_j}d^{-\frac 12}).$$
Hence $\Psi_{\mathcal{G}^{\varphi}}(J(f))$ has the matrix representation:
$$\Phi(\Psi_{\mathcal{G}^{\varphi}}(J(f)))=diag((1_{U_1}d^{-\frac 12})\Psi_{\mathcal{G}}(fd^{\frac 12}),\cdots, (1_{U_i}d^{-\frac 12})\Psi_{\mathcal{G}}(fd^{\frac 12}),\cdots, (1_{U_N}d^{-\frac 12})\Psi_{\mathcal{G}}(fd^{\frac 12})).$$

		For $u\in U_i$, the range-fibre sum of $J(f)$ is
		\[
		\sum_j\sum_{\gamma\in \mathcal{G}^u}\frac{|f(\gamma)|1_{U_j}(s(\gamma))}{\sqrt{d(u)d(s(\gamma))}}
		=\frac{1}{\sqrt{d(u)}}\sum_{\gamma\in \mathcal{G}^u}|f(\gamma)|\sqrt{d(s(\gamma))}
		\le\sqrt{N}\sum_{\gamma\in \mathcal{G}^u}|f(\gamma)|.
		\]
		The source-fibre estimate is identical. This proves $\|J(f)\|_{I,\mathcal{G}^\varphi}\le \sqrt{N}\|f\|_{I,\mathcal{G}}.$
		
	Note that $W_{\rho}\xi=\oplus_{i=1}^N(\rho(1_{U_i} d^{-\frac 12})\xi)$ and $W_{\rho}^{\ast}(\oplus_{i=1}^N\xi_i)=\sum\limits_{i=1}^N(\rho(d^{-\frac 12})\xi_i)$ for each $\xi\in H$  and $\oplus_{i=1}^N\xi_i\in \mathcal{L}$.
Thus equation (\ref{eq4}) follows from direct matrix multiplications. \end{proof}

	In order to modify the equation (\ref{eq4}), we have the following result. For $F\in C_c(\mathcal{G}^\varphi)$, write
		\[
		D_F:=\pi_\rho(F)-\pi_\rho(\Psi_{\mathcal{G}^\varphi}(F)).
		\]
	
	\begin{lemma}\label{Le10}
		For $f\in C_c(\mathcal{G})$, let $\kappa_f=fd^{-\frac 12}\in C_c(\mathcal{G})$.  Then $\|\kappa_f\|_{I,\mathcal{G}}\le\|f\|_{I,\mathcal{G}}$, and
		\begin{equation}
			W_\rho^*D_{J(\kappa_f)}W_\rho
			=(\rho(f)-\rho(\Psi_{\mathcal{G}}(f)))\rho(d^{-\frac 12}).
			\label{eq6}
		\end{equation}
		Moreover, the operator $D_{J(\kappa_f)}$
				commutes with $\pi_\rho(e_\varphi)$.
	\end{lemma}
	
	\begin{proof} The inequality for $\|\cdot\|_{I}$ follows from the definitions. Since $W_{\rho}$ is an isometry, equation \eqref{eq6} follows from  \eqref{eq4}.

Noting that $\pi_\rho(e_\varphi)$ is the range projection of $W_{\rho}$,  we also obtain that $\pi_\rho(J(\kappa_f))$ commutes with $\pi_\rho(e_\varphi)$ from the first equality in equation \eqref{eq4}. Given $g\in C_c(\mathcal{G})$, by calculation, we have $$\Psi_{\mathcal{G}^{\varphi}}(J(g))(i,u,i)=d^{-\frac 12}(u)\sum\limits_{\gamma\in \mathcal{G}^{u}}f(\gamma)d^{\frac 12}(s(\gamma))$$ for each $(i,u,i)\in (\mathcal{G}^{\varphi})^{(0)}$, which can implies that $e_{\varphi}\Psi_{\mathcal{G}^{\varphi}}(J(g))=\Psi_{\mathcal{G}^{\varphi}}(J(g))e_{\varphi}$. Hence $D_{J(\kappa_f)}$ commutes with $\pi_\rho(e_\varphi)$.
	\end{proof}

The following lemma shows that if $\rho$ is a (left) regular representation of $C_c(\mathcal{G})$, then $\pi_{\rho}$ is unitarily equivalent to a regular representation of $C_c(\mathcal{G}^{\varphi})$.
	
	 \begin{lemma}\label{Le11}
			Let $x\in X$ and put $u=\varphi(x)$. If $\rho=\pi_u: C_c(\mathcal{G})\rightarrow B(l^2(\mathcal{G}_u))$ is the regular representation at $u$, then $\pi_{\rho}$ is  unitarily equivalent to the regular representation $\pi_x$ of $C_c(\mathcal{G}^{\varphi})$ at $x$.

under the sheet decomposition of $\varphi$,
			there is a canonical unitary
			\[
			\ell^2\bigl((\mathcal{G}^\varphi)_x\bigr)\longrightarrow
			\bigoplus_{i=1}^N \pi_u(1_{U_i})\ell^2(\mathcal{G}_u),
			\]
		which identifies the regular representation of $\mathcal{G}^\varphi$ at $x$
			with   $\pi_{\pi_{u}}$.
			
	\end{lemma}
	
	\begin{proof} By the decomposition of $\varphi$, there exists a unique $j_0$ ($1\leq j_0\leq N$) such that $x\in X_{j_0}$, so $u\in U_{j_0}$ and $x=\sigma_{j_0}(u)$. From the identification of $\mathcal{G}_{\mathcal{U}}$ with $\mathcal{G}^{\varphi}$, we identify $x$ with $(j_0,u,j_0)$. Hence $(\mathcal{G}^{\varphi})_x=\{(i,\gamma,j_0): \gamma\in \mathcal{G}^{U_i}_u\}$.

Let $H=l^2(\mathcal{G}_u)$ and let $\{\delta_{\gamma}: \gamma\in \mathcal{G}_u\}$ be the canonical orthogonal basis of $H$. For each $i$, we check that $\rho(1_{U_i})\delta_{\gamma}=\delta_{\gamma}$ for $r(\gamma)\in U_i$, and $\rho(1_{U_i})\delta_{\gamma}=0$ for $r(\gamma)\notin U_i$. Hence $\{\delta_{\gamma}: \gamma\in \mathcal{G}^{U_i}_u\}$ is the canonical orthogonal basis of $\rho(1_{U_i})H$. Thus $\mathcal{L}:=\oplus_{i=1}^N(\rho(1_{U_i})H)$ has the canonical orthogonal basis $\{(0,\cdots, 0, \delta_{\gamma_i}, 0,\cdots, 0)^{tr}\in \mathcal{L}: \gamma_i\in \mathcal{G}^{U_i}_u, 1\leq i\leq N\}$, where $(0,\cdots, 0,\delta_{\gamma_i}, 0,\cdots, 0)^{tr}$ denotes the column vector.

		Define the map $U: l^2((\mathcal{G}^{\varphi})_x)\rightarrow \mathcal{L}$ by
		$$
		\delta_{(i,\gamma,j_0)}\longmapsto (0,\cdots, 0,\delta_{\gamma},0,\cdots,0)^{tr}.$$
		 We can check that $U$ is a unitary operator and satisfies that $U^*\pi_{\rho}(\Phi(F))U=\pi_x(F)$ for each $F\in C_c(\mathcal{G}^{\varphi})$. That is, $U$ is a unitary intertwining $\pi_{\rho}$ and the regular representation $\pi_x$ of $C_c(\mathcal{G}^{\varphi})$ at $x$ under the identification in Lemma \ref{le0}.
	\end{proof}
	
	In the rest of this section, we compare Property (T) for an ample groupoid $\mathcal{G}$ and its ampliation $\mathcal{G}^{\varphi}$ along a surjective local homeomorphism $\varphi: X\to \mathcal{G}^{(0)}$.
We adopt the notations as before and have the following main result.

\begin{theorem}\label{Th4} Let $\mathcal{G}$ be an ample groupoid, $X$ be a compact Hausdorff space and let $\mathcal{G}^{\varphi}$ be the the ampliation of $\mathcal{G}$ with respective to  a surjective local homeomorphism $\varphi$ from $X$ onto $\mathcal{G}^{(0)}$. Then  $\mathcal{G}$ has property $(T)$ if and only if $\mathcal{G}^{\varphi}$ has  property $(T)$.

Moreover, $(\mathcal{G},\mathcal{R}_\mathcal{G})$ has Property {\rm(T)} if and only if $(\mathcal{G}^\varphi,\mathcal{R}_{\mathcal{G}^\varphi})$ has Property {\rm(T)}.

\end{theorem}

This theorem consists the following two propositions.
	
	\begin{proposition}\label{Pro3}  Suppose $(K,c)$ is a Kazhdan pair for $\mathcal{G}$ with respect to all representations, then $(\widetilde{L}, c_{\varphi})$ is a  Kazhdan pair for $\mathcal{G}^{\varphi}$ with respect to all representations, where $\widetilde{L}=q^{-1}(K\cup \mathcal{G}^{(0)})$ and
		$c_{\varphi}=\frac{1}{\sqrt2}\min\left(1,\frac{c}{N}\right)$.

		The same assertion holds for the family of regular representations.
	\end{proposition}
	
	\begin{proof}  We can assume that $K$ contains $\mathcal{G}^{(0)}$. For a representation $\pi$ of $C_c(\mathcal{G}^{\varphi})$,  by Lemma~\ref{Le7}, we can assume that $\pi$ has the form $\pi_\rho$ on $\mathcal L$ which is induced by a representation $\rho$ of $C_c(\mathcal{G})$ on $H$.
Given $\xi\in\mathcal L_{\pi_\rho}=(\mathcal L^{\pi_\rho})^{\perp}$, let $\xi=\xi_1+\xi_2$ for $\xi_1=\pi_\rho(e_\varphi)\xi$ and $\xi_2=(1-\pi_\rho(e_\varphi))\xi$.
		Since $\mathcal L^{\pi_\rho}=S_{\rho}(H^{\rho})\subseteq S_{\rho}(H)=Ran(\pi_\rho(e_\varphi))$, we have $\xi_i\in \mathcal L_{\pi_\rho}$ for $i=1,2$. We have two cases.
		
		Case 1: $\|\xi_2\|\geq\frac{1}{\sqrt2}\|\xi\|$. Note that we have obtained that $\|e_\varphi\|_{I,\mathcal{G}^\varphi}=1$ and $\operatorname{supp}(e_\varphi)\subseteq \widetilde{\mathcal{G}^{(0)}}\subseteq \widetilde{K}$. Also since $\Psi_{\mathcal{G}^\varphi}(e_\varphi)=1_{X}$, we have
		$$
		\|\pi_\rho(e_\varphi)\xi-\pi_\rho(\Psi_{\mathcal{G}^\varphi}(e_\varphi))\xi\|=\|\xi_1-\xi\|=\|-\xi_2\|\geq \frac{1}{\sqrt2}\|\xi\|\geq  c_{\varphi}\|\xi\|.$$
				
		Case 2: $\|\xi_1\|\ge\frac{1}{\sqrt2}\|\xi\|$. Let $\eta=W_\rho^*\xi_1$. Since $W_{\rho}(\rho(d^{\frac 12})H^{\rho})=S_{\rho}(H^{\rho})=\mathcal L^{\pi_\rho}$, we have that $\eta\perp\rho(d^{\frac 12})(H^\rho)$.
Put $\zeta=\rho(d^{-\frac12})\eta$ and let $\zeta_0$ be the orthogonal projection of $\zeta$ onto $H_\rho$. Thus $(\zeta-\zeta_0)\in H^\rho$. In this case, we have $$\|\rho(d^{\frac 12})\zeta_0\|^2=\|\eta\|^2+\|\rho(d^{\frac12})(\zeta-\zeta_0)\|^2\geq \|\eta\|^2=\|\xi_1\|^2, $$
which implies that
	\[
		\|\zeta_0\|\ge\frac{\|\xi_1\|}{\sqrt N}.
		\]
		Applying the Kazhdan condition for $\mathcal{G}$ to $\zeta_0$, we obtain $f\in C_c(\mathcal{G})$ with
		\[
		\operatorname{supp}(f)\subseteq K,\qquad \|f\|_{I,\mathcal{G}}\le 1,
		\]
		such that
		$$
		\|(\rho(f)-\rho(\Psi_\mathcal{G}(f)))(\zeta_0)\|
		\ge c\|\zeta_0\|\ge\frac{c}{\sqrt N}\|\xi_1\|.	
		$$
		Since the left-hand side is unchanged when $\zeta_0$ is replaced by $\zeta$, we have
		\[
		\|(\rho(f)-\rho(\Psi_{\mathcal{G}}(f)))\zeta\|\ge\frac{c}{\sqrt N}\|\xi_1\|.
		\]

		Let $\kappa_f$ be defined by Lemma~\ref{Le10}. By Lemmas~\ref{Le9} and \ref{Le10}, we have $W_{\rho}^*D_{J(\kappa_f)}\xi_1=(\rho(f)-\rho(\Psi_\mathcal{G}(f)))\zeta$, thus
		\begin{equation}
			\|D_{J(\kappa_f)}\xi_1\|\ge\frac{c}{\sqrt N}\|\xi_1\|.
			\label{eq9}
		\end{equation}
		Since the operator $D_{J(\kappa_f)}$ commutes with $\pi_\rho(e_\varphi)$, we have that $D_{J(\kappa_f)}\xi_1$ and $D_{J(\kappa_f)}\xi_2$ are orthogonal, so
		$$
		\|D_{J(\kappa_f)}\xi\|\ge\|D_{J(\kappa_f)}\xi_1\|\ge\frac{c}{\sqrt N}\|\xi_1\|\geq\frac{c}{\sqrt{2N}}\|\xi\|.$$
		Note that $
		\|J(\kappa_f)\|_{I,\mathcal{G}^{\varphi}}\le\sqrt N\|\kappa_f\|_{I, \mathcal{G}}\le\sqrt N\|f\|_{I, \mathcal{G}}\leq \sqrt N$.
		Let
		$$
		F:=\frac{1}{\sqrt N}J(\kappa_f).$$
Then $F\in C_c(\mathcal{G}^{\varphi})$
		has $I$-norm at most one, is supported in $\tilde K$, and satisfies
		\[
		\|D_F\xi\|\ge\frac{c}{\sqrt2 N}\|\xi\|\geq c_{\varphi}\|\xi\|.\]
		Combining the two cases, we show that $(\widetilde{K}, c_{\varphi})$ is a Kazhdan pair for $\mathcal{G}^{\varphi}$.
		
	 	For regular representations, the same argument applies when the regular representation of $\mathcal{G}^\varphi$ is identified with the appropriate $\pi_\rho$.\end{proof}

	\begin{proposition}\label{Pro4}
		Suppose $(L,\delta)$ is a Kazhdan pair for $\mathcal{G}^{\varphi}$ with respect to all representations. Then $(q(L),\frac{\delta}{N^{3/2}})$ is a Kazhdan pair for $\mathcal{G}$ with respect to all representations.

		The same assertion holds for regular representations.
	\end{proposition}
	
	\begin{proof}
		Fix a representation $\rho: C_c(\mathcal{G})\to B(H)$ and a non-zero vector $\xi\in  H_\rho$. We consider the induced representation $\pi_\rho: C_c( \mathcal{G}^{\varphi})\to B(\mathcal L)$. Let $Q$ be the orthogonal projection of $\mathcal{L}$ onto
		$
		\mathcal L^{\pi_\rho}\,\,\,(=S_\rho(  H^\rho)),
		$
		and define
		$
		\eta:=(1-Q)S_\rho\xi.
		$
		Then $\eta\in\mathcal L_{\pi_\rho}$. Since
		$
		S_\rho^*S_\rho=\rho(d)\ge 1,
		$
		we have
		\[
		\|\eta\|
		=\inf_{h\in  H^\rho}\|S_\rho(\xi)-S_{\rho}(h)\|
		\ge\inf_{h\in  H^\rho}\|\xi-h\|=\|\xi\|.
		\]

		Since $(L,\delta)$ is a Kazhdan pair for $\mathcal{G}^{\varphi}$, for the vector $\eta$,  we choose $F\in C_c(\mathcal{G}^{\varphi})$ such that $\operatorname{supp}(F)\subseteq L$, $\|F\|_{I,\mathcal{G}^{\varphi}}\le 1$, and $$
		\|D_F\eta\|\ge\delta\|\eta\|,$$
		where $D_F:=\pi_\rho(F)-\pi_\rho(\Psi_{\mathcal{G}^{\varphi}}(F))$.  Since $D_F$ vanishes on the constant-vector space, we have
		$
		D_F\eta=D_F S_\rho\xi.
		$

		In the matrix representation, let $F=(f_{ij})$. For each $i$, define
		$$
		g_i:=\sum_{j=1}^N f_{ij}\in C_c(\mathcal{G}).
		$$
		Using the matrix-corner description of Proposition~\ref{Pro2}, one computes
		$$
		D_F S_\rho\xi=\oplus_{i=1}^n (\rho(g_i)-\rho(\Psi_\mathcal{G}(g_i))\bigr)\xi.	$$
		Hence,
		\[
		\sum_{i=1}^N\bigl\|(\rho(g_i)-\rho(\Psi_\mathcal{G}(g_i)))\xi\bigr\|^2=\|D_F S_\rho\xi\|^2=\|D_F \eta\|^2\ge\delta^2\|\xi\|^2.
		\]
		Thus for some $i$,
		\begin{equation}
			\bigl\|(\rho(g_i)-\rho(\Psi_\mathcal{G}(g_i)))\xi\bigr\|\ge\frac{\delta}{\sqrt N}\|\xi\|.\label{eq10}
		\end{equation}
		Clearly,
		$
		\operatorname{supp}(g_i)\subseteq q(L).
		$
		The range-fibre part of $\|g_i\|_{I,\mathcal{G}}$ is bounded by one using the row estimate in Proposition~\ref{Pro2}. For the source-fibre part,
		\[
		\sup_u\sum_{h\in \mathcal{G}_u}|g_i(h)|
		\le\sum_{j=1}^N\sup_u\sum_{h\in \mathcal{G}_u}|f_{ij}(h)|
		\le N.
		\]
		Thus
		$
		\|g_i\|_{I,\mathcal{G}}\le N.
		$
		Replacing $g_i$ by $g_i/N$ in equation \eqref{eq10} gives the asserted Kazhdan constant $\delta/N^{3/2}$.
		\end{proof}
		
	%

\begin{corollary} \label{Cor2}
		Let $\mathcal{G}$ be ample with compact unit space and let $U\subseteq \mathcal{G}^{(0)}$ be an arbitrary full and clopen subset. Then
		$\mathcal{G}$ has topological Property {\rm(T)} if and only if $\mathcal{G}|_U$  has topological Property {\rm(T)}.
		
		Likewise, $(\mathcal{G},\mathcal{R}_\mathcal{G})$  has Property {\rm(T)} if and only if $(\mathcal{G}|_U,\mathcal{R}_{\mathcal{G}|_U})$ has Property {\rm(T)}.
	\end{corollary}
	
	\begin{proof}  Let $\mathcal{H}=\mathcal{G}|_U$. Then it is an ample \'{e}tale groupoid with compact unit space $U$. From Theorem \ref{Th4}, this corollary follows from the claim that $\mathcal{G}$ is isomorphic to some ampliation $\mathcal{H}^{\varphi}$ of $\mathcal{H}$.

		For each $x\in \mathcal{G}^{(0)}\setminus U$, it follows from the fullness of $U$ that there exists $\gamma_x\in \mathcal{G}$ with $r(\gamma_x)=x$, $ s(\gamma_x)\in U.$
		Since $\mathcal{G}$ is ample, choose a compact open bisection neighbourhood $B_x$ of $\gamma_x$ in $\mathcal{G}$  such that $s(B_x)\subseteq U$, $r(B_x)\subseteq \mathcal{G}^{(0)}\setminus U$. Thus $\{r(B_{x}):x\in \mathcal{G}^{(0)}\setminus U\}$ is a clopen cover of $\mathcal{G}^{(0)}\setminus U$.
		By compactness and zero dimensionality of $\mathcal{G}^{(0)}\setminus U$, there is a finite pairwise disjoint clopen partition
		\[
		\mathcal{G}^{(0)}\setminus U=C_1\sqcup\cdots\sqcup C_m
		\]
		such that each $C_i$ is contained in some $r(B_x)$. Restricting the bisection $B_x$ to $r^{-1}(C_i)$, we may assume
		$
		r(B_i)=C_i$, $s(B_i)\subseteq U.
		$
		Set also
		\[
		C_0:=U,\qquad B_0:=U.
		\]
		We have a clopen partition $\{C_i: i=0,1,2,\cdots,m\}$ of $\mathcal{G}^{(0)}$. Let $\beta_i: C_i\to B_i$ be the inverse of the range map $r|_{B_i}: B_i\to C_i$ for $i=0,1,2,\cdots,m$.  Consequently, we have a continuous map $\beta: \mathcal{ G}^{(0)}\rightarrow \mathcal{G}$ such that $\beta|_{C_i}=\beta_i$ for each $i$. Clearly, $\beta|_{U}=\beta_0=id_U$.

		Define $\varphi=s\beta: \mathcal{G}^{(0)}\to \mathcal{G}^{(0)}$. Then it is a  surjective local homeomorphism from $\mathcal{G}^{(0)}$ onto $U$ such that $\varphi|_U=id_U$. Define
		$
		\Theta:\mathcal{G}\to \mathcal{H}^{\varphi}
		$
		by
		$$
		\Theta(\gamma)=\bigl(r(\gamma),\beta(r(\gamma))^{-1}\gamma\beta(s(\gamma)),s(\gamma)\bigr). $$
		Then  $\Theta$ is a continuous groupoid isomorphism whose inverse is defined by
		\[
		(x,k,y)\mapsto \beta(x)k\beta(y)^{-1}.
		\]
		Hence $\mathcal{G}$ is isomorphic to the amplication $\mathcal{H}^{\varphi}$. We establish the claim.
\end{proof}



From the definition of Morita equivalence and Theorem \ref{Th4}, we have the following main result.

			\begin{theorem} \label{Th5}

				Let $\mathcal{G}$ and $\mathcal{H}$ be ample \'{e}tale groupoids with compact unit spaces. If $\mathcal{G}$ and $\mathcal{H}$ are Morita equivalent, then
$\mathcal{G}$ has topological Property {\rm(T)} if and only if $\mathcal{H}$  has topological Property {\rm(T)}.
		
		Likewise, $(\mathcal{G},\mathcal{R}_\mathcal{G})$  has Property {\rm(T)} if and only if $(\mathcal{H},\mathcal{R}_{\mathcal{H}})$ has Property {\rm(T)}. \end{theorem}

			\begin{remark} Combining the Farsi--Kumjian--Pask--Sims equivalence theorem \cite[Theorem 3.12]{Farsi} with Theorem~\ref{Th5}, for ample Hausdorff \'etale groupoids with compact unit spaces, we have that  topological Property {\rm(T)} is invariant under each of the following equivalent notions:
				similarity, Morita equivalence, Renault equivalence, isomorphic ampliations, stable isomorphism, Kakutani equivalence, weak Kakutani equivalence. The same statement holds for Property {\rm(T)} with respect to the family of regular representations.
			\end{remark}

%
%
%
%

		\end{document}